\documentclass[oneside, 12pt]{amsart}
\usepackage{amscd, amssymb, amsmath, mathrsfs}
\usepackage[english]{babel}
\usepackage{booktabs}
\usepackage{tikz-cd}
\usepackage[pdftex, colorlinks=true,  citecolor=blue, linkcolor=blue, linktocpage=true]{hyperref}

\newcommand\mylabel[1]{\label{#1}\marginpar{\vspace{-1ex}\medskip\medskip\footnotesize \tt #1}}
\renewcommand\mylabel[1]{\label{#1}}
\newcommand{\mydate}{
\number\day\space
\ifcase\month \or January\or February\or March\or April\or May\or June\or July\or August\or September\or October\or November\or December\fi 
\space\number\year}

\newtheorem{theorem}{Theorem}[section]
\newtheorem*{maintheorem}{Theorem}
\newtheorem{lemma}[theorem]{Lemma}
\newtheorem{proposition}[theorem]{Proposition}
\newtheorem{corollary}[theorem]{Corollary}

\theoremstyle{definition}

\newtheorem*{acknowledgement}{Acknowledgement}

\theoremstyle{remark}

\newcommand{\ZZ}{\mathbb{Z}}
\newcommand{\QQ}{\mathbb{Q}}

\newcommand{\FF}{\mathbb{F}}

\newcommand{\PP}{\mathbb{P}}
\renewcommand{\AA}{\mathbb{A}}
\newcommand{\GG}{\mathbb{G}}

\newcommand{\ideala}{\mathfrak{a}}

\newcommand{\shM}{\mathscr{M}}

\newcommand{\shL}{\mathscr{L}}

\newcommand{\shMbar}{\bar{\shM}}

\newcommand{\alg}{\text{\rm alg}}

\newcommand{\Aut}{\operatorname{Aut}}

\newcommand{\Frac}{\operatorname{Frac}}

\newcommand{\Gal}{\operatorname{Gal}}

\newcommand{\I}{\text{\rm I}}

\newcommand{\lra}{\longrightarrow}

\newcommand{\maxid}{\mathfrak{m}}

\renewcommand{\O}{\mathscr{O}}

\newcommand{\Pic}{\operatorname{Pic}}

\newcommand{\Proj}{\operatorname{Proj}}

\newcommand{\quadand}{\quad\text{and}\quad}

\newcommand{\ra}{\rightarrow}

\newcommand{\Reg}{\operatorname{Reg}}

\newcommand{\Sing}{\operatorname{Sing}}

\newcommand{\Spec}{\operatorname{Spec}}

\newcommand{\val}{\operatorname{val}}

\begin{document}

\title[Elliptic curves over the rational numbers]
      {Elliptic curves over the rational numbers with semi-abelian reduction and two-division points}.

\author[Stefan Schr\"oer]{Stefan Schr\"oer}
\address{Mathematisches Institut, Heinrich-Heine-Universit\"at,
40204 D\"usseldorf, Germany}
\curraddr{}
\email{schroeer@math.uni-duesseldorf.de}


\dedicatory{Second revised Version, 12 October 2020}

\begin{abstract}
We classify elliptic curves over the   rationals
whose N\'eron model over the integers is semi-abelian,
with good reduction at $p=2$, and whose Mordell--Weil group
contains an element of order two that stays non-trivial at $p=2$.
Furthermore, we describe those curves where the element of order two is narrow, or
where another element of order two exists, and also express our findings
in terms of Deligne--Mumford stacks of   pointed curves of genus one.
\end{abstract}

\maketitle
\tableofcontents

\section*{Introduction}
\mylabel{introduction}

Fontaine \cite{Fontaine 1985} proved that there are no abelian varieties $A_\QQ\neq 0$  over the field of rationals
that have good reduction everywhere. The   case of elliptic curves  $E_\QQ$
is due to Ogg \cite{Ogg 1966}, who attributed the result to Tate.  His proof proceeds by showing that any  Weierstra\ss{} equation 
$
y^2+a_1xy+a_3y=x^3+a_2x^2+a_4x+a_6
$
with integral coefficients has discriminant $\Delta\neq \pm 1$. Ogg actually classified, up to isomorphism,
those $E_\QQ$ where bad reduction occurs only at $p=2$. It turns out that there are twelve cases,
all of the form $y^2=x^3+mx^2+nx$ with discriminant $\Delta=\pm 2^v$,
for certain coefficients $-12\leq m,n\leq 12$ and   exponents  $6\leq v\leq 15$. 
One observes that there is always a 2-division point, that is, a non-zero rational point of order two, with coordinates $x=y=0$. 

Similar questions where treated in many subsequent papers. We just mention a few:
Ogg himself analyzed elliptic curves with
discriminant $\pm 2^\nu\cdot 3$ and $\pm 2^\nu\cdot 3^2$ in  \cite{Ogg 1967}.
The case $\Delta=\pm 2^\nu\cdot p$ was studied  by Ivorra \cite{Ivorra 2004},   those
with $\Delta=p_1p_2$   by Sadek   \cite{Sadek 2014}.
Ogg's results where extended to imaginary quadratic number fields by Setzer \cite{Setzer 1978},
Stroecker \cite{Stroecker 1983} and Kida \cite{Kida 2001}.
Zhao \cite{Zhao 2013} studied elliptic curves $E_K$ with good reduction everywhere over real quadratic number fields,
and Takeshi \cite{Takeshi 2015} investigated the case of cubic number fields.
The reduction behavior of an elliptic curve $E_\QQ$ with a 2-division point was studied by
Hadano \cite{Hadano 1974}.

The goal of this paper is to classify elliptic curves $E_\QQ$
that have a 2-division point $P\in E(\ZZ)$, and  satisfy certain additional conditions. 
Here $E\ra\Spec(\ZZ)$ is the N\'eron model, and $E(\ZZ)=E(\QQ)=E_\QQ(\QQ)$
is the \emph{Mordell--Weil group}, which is a finitely generated abelian group.
My personal motivation is to study the arithmetic
of families of elliptic surfaces, in particular Enriques surfaces \cite{Schroeer 2020}, but the results
seem to be of independent interest.
The additional conditions considered here are as follows: \emph{We stipulate that the relative group scheme $E$
is semi-abelian, that good reduction occurs at $p=2$, and that the $2$-division point 
$P\in E(\ZZ)$  stays non-trivial in $E(\FF_2)$.}

It turns out that there are infinitely many isomorphism classes. We show    
that they correspond to pairs of integers $(m,n)$ with $n$ odd and $\gcd(4m+1,n)=1$, via the Weierstra\ss{} equations
$y^2+xy=x^3+mx^2+nx$ (see Theorem \ref{equation necessary}).
We then describes the situation when there is further 2-division, which is our first main result:

\begin{maintheorem}
{\rm (See Theorem \ref{additional 2-division})}
The following are equivalent:
\begin{enumerate}
\item There is another element $Q\neq P$ of order two in $E(\ZZ)$.
\item The relative group scheme $E[2]$ is isomorphic to $\ZZ/2\ZZ\oplus\mu_2$.
\item For every prime $p>0$, the fiber $E\otimes\FF_p$ has   even  Kodaira symbol.
\item We have $n = -d(4m +1+16d)$ for some odd   $d$ with $\gcd( 4m+1,d)=1$.
\end{enumerate}
\end{maintheorem}

I find it quit remarkable that the existence of $Q\neq P$ can be seen as  a condition on   relative group schemes,
or   Kodaira symbols, or the arithmetic of the numbers $m$ and $n$.

Now recall that $P\subset E$ is called \emph{narrow} it it passes through the same  irreducible components
of the closed fibers as the zero-section $O\subset E$, compare for example \cite{Schuett; Shioda 2019},  Section 6.7.
If we  stipulate that our $P\subset E$ is  narrow,   the  above list becomes finite.
The second main results of the paper is:

\begin{maintheorem}
{\rm (See Theorem \ref{classification order two})}
Up to isomorphism, there  are exactly two elliptic curves $E_\QQ$ such that its N\'eron model $E$ has
the following properties:
\begin{enumerate}
\item The structure morphism $E\ra\Spec(\ZZ)$ is semi-abelian.
\item The closed fiber $E\otimes\FF_2$ is an   elliptic curve.
\item There is a narrow element $P\in E(\ZZ)$ of order two whose image  in
$E(\FF_2)$ is non-zero, and another element $Q\neq P$ of order two.
\end{enumerate}
These  curves are given by $y^2+xy = x^3+4nx^2+nx$ for the coefficient   $n=\pm1$.
\end{maintheorem}

Taking the quotient by the subgroup scheme $\mu_2\subset E$  yields two other elliptic
curves, which must be of the form $y^2+xy=x^3+m'x^2+n'x$.
Somewhat surprisingly, they can be characterized via elements $R\in E(\ZZ)$ of order four,
which is the content of the third main result (see Theorem \ref{classification order four}).
We also reformulate some of our findings in terms of of $\ZZ$-valued objects
in the Deligne--Mumford stack $\shMbar_{g,r}$ of stable pointed curves of genus $g=1$.

\medskip
The paper is organized as follows:
Section \ref{Preliminaries} contains basic facts on N\'eron models of elliptic curves,
and elementary computations with Weierstra\ss{} equations.
In Section \ref{Geometry to equations}, we show how certain geometric assumptions on N\'eron models
lead to our Weierstra\ss{} equations. A detailed analysis of the resulting Weierstra\ss{} models
is given in Section \ref{Analysis}, which also contains  the result on the existence of an additional element of order two.
In Section \ref{Quotients} we apply Velu's formula to form some quotients by 2-division points.
Section \ref{Classification} contains our classification results that use additional narrowness assumptions.
The final Section \ref{Reinterpretation} gives reinterpretations in terms of Deligne--Mumford stacks.

\begin{acknowledgement}
I wish to thank the referee for valuable suggestions.
This research was conducted in the framework of the   research training group
\emph{GRK 2240: Algebro-geometric Methods in Algebra, Arithmetic and Topology}, which is funded
by the Deutsche Forschungsgemeinschaft. 
\end{acknowledgement}

\section{Preliminaries}
\mylabel{Preliminaries}

Here we discuss basic   facts  on N\'eron models of elliptic curves and 
make some useful elementary computations with Weierstra\ss{} equations. 
We refer to the monograph of Bosch, L\"utkebohmert and Raynaud  for more details
(\cite{Bosch; Luetkebohmert; Raynaud 1990},  Section 1.5 in particular).
Let $R$ be a Dedekind ring, which we regard as a ground ring, $F=\Frac(R)$ its field of fractions,
and
\begin{equation}
\label{weq}
E_F:\quad y^2+a_1xy + a_3y=x^3+a_2x^2+a_4x+a_6
\end{equation}
be an elliptic curve with coefficients $a_i\in F$.
The homogenization of the Weierstra\ss{} equation yields a closed embedding
$E_F\subset\PP^2_F$, where the origin  $O_F\in E_F$ corresponds to the point $(0:1:0)\in\PP^2(F)$.

Let $X'$ be the resolution of singularities for the closure  of the subset $E_F\subset \PP^2_R$, which we assume
to exist. Note that almost all closed fibers are irreducible. 
Let $X'\ra X$ be a relatively minimal model, which arises by a succession of contractions of $(-1)$-curves.
The closure of the origin defines a section $O\subset X$. 
According to \cite{Schroeer 2001}, Corollary 2.3 
there is  a contraction $X\ra Y$ that contracts the irreducible components of the closed fibers that
are disjoint from $O$. Then the image  of $O\subset X$ on $Y$ is an relatively ample effective Cartier divisor, and it follows that 
$$
Y=P(X,O)=\Proj \bigoplus_{t\geq 0} H^0(X,\O_{X}(tO)).
$$
Note that $X\ra Y$ is the minimal resolution of singularities.
The schemes 
$X$ and $Y$ are called the   \emph{minimal   model} and the \emph{Weierstra\ss{} model}
of the elliptic curve $E_F$, respectively. Indeed, the scheme $X$, and hence also $Y$, neither  depend on the chosen
Weierstra\ss{} equation  nor the resolution $X'$, up to isomorphism: If $Z_1\subset \PP^2_R$ is the original closure,
and $Z_2\subset\PP^2_R$ is the closure of $E_F$ with respect to another Weierstra\ss{} equation, we also find
some possibly singular compactification $X''$ of $E_F$ that dominates $X'$. However, by   \cite{Stacks Project}, Lemma 0AHH
we   may  assume that $X''$ is regular, and that $X''\ra X'$ is a succession of contractions of $(-1)$-curves.
It thus remains to verify that any two regular compactifications
$X',X''$ of $E_F$ give the same $X$.  This holds by \cite{Lichtenbaum 1968}, Theorem 4.4, because $h^1(\O_{E_F})\neq 0$.

Each  fiber  of $Y\ra\Spec(R)$ is either an  elliptic curves, or a twisted from of 
the rational nodal curve  or the rational cuspidal curve.
The \emph{locus of non-smoothness} $\Sing(Y/R)$, which is defined by the first Fitting ideal
for the sheaf $\Omega^1_{Y/R}$ (compare \cite{Fanelli; Schroeer 2018}, Section 2), comprises   finitely many closed points.
The locus of smoothness $E=\Reg(X/R)$ in the minimal   model is  the \emph{N\'eron model} of the
elliptic curve $E_F$. If the residue fields $R/\maxid$ are perfect, this follows from \cite{Bosch; Luetkebohmert; Raynaud 1990}, Section 1 in Section 1.5,
together with Theorem 1 in Section 7.2. The general case follows from \cite{Szydlo  2004}, Proposition 7.1.1.
The N\'eron model has the  property that the restriction map
$E(R)\ra E(F)$ is bijective. The resulting abelian group
$$
E(R) = E(F)=E_F(F) 
$$
is called the \emph{Mordell--Weil group}.
The N\'eron model $E\ra\Spec(R)$ acquires in a canonical way the structure of  a relative group scheme.
We write $O\subset E$ for the zero-section.
For each point $a\in \Spec(R)$, the fiber $E_a=E\otimes\kappa(a)$ is an algebraic group scheme.
We write $E_a^0\subset E_a$ for the connected component of the origin,
and $\Phi_a=E_a/E_a^0$ for the resulting   group scheme of components.
If  the fiber $E_a$ is an elliptic curve, one says that $E_F$ has \emph{good reduction} at $a\in\Spec(R)$.
If $E_a^0$ is a twisted form of $\GG_m\otimes\kappa(a)$, one says that $E_F$ has \emph{multiplicative reduction}.
One says that the N\'eron model $E$ is \emph{semi-abelian} if the elliptic curve
$E_F$ has good reduction or multiplicative reduction at all closed points.

An element $P\in E(R)$ is called \emph{narrow} if its images in $\Phi_a$ vanishes,
for all closed points $a\in \Spec(R)$. In other words, the section $P,O\subset E$ pass through
the same irreducible components inside each closed fiber $E_a$.

There is an open covering $ D(f_1)\cup\ldots\cup D(f_n)$ of $\Spec(R)$ such that 
each restriction
$Y\otimes R[1/f_i]$ is   defined by some   equation \eqref{weq} with coefficients $a_i\in R[1/f_i]$
as a family of cubic curves in $\PP^2\otimes R[1/f_i]$.
In general, it is impossible to find a global  embedding $Y\subset\PP^2$.
However, making  a substitution $x=u^2x'$ and $y=u^3y'$ with $u\in F^\times$ we obtain a Weierstra\ss{} equation \eqref{weq} with   coefficients $a_i\in R$. Such a Weierstra\ss{} equation is called \emph{globally minimal} if the resulting family of cubics
is isomorphic to $Y\ra\Spec(R)$. 

If the Dedekind ring $R$ is local, that is, a discrete valuation ring, 
the \emph{Tate Algorithm} \cite{Tate 1972a} produces a minimal Weierstra\ss{} equation, at least if the residue field is perfect.
Every round of the algorithm consists of twelve steps, each involving a  
change   of coordinates $x=u^2x'+r$ and $y=u^3y'+su^2x'+t$, with $u\in R^\times$
and $r,s,t\in R$, to produces new equations. If the situation  $\val(a_i)\geq i$ arises,
one replaces the coefficients by $a_i/\pi^i$, where $\pi\in R$ is the uniformizer,
and starts a new round. If not, the Weierstra\ss{} equation is minimal and indeed describes $Y$.
It the ring $R$ is factorial, one actually may run the Tate Algorithm locally and thus obtains
a globally minimal Weierstra\ss{} equation.

\begin{lemma}
\mylabel{globally minimal}
Suppose the Weierstra\ss{} equation \eqref{weq} has coefficients $a_i\in R$ and that the values $c_4,\Delta\in R$ generate the unit ideal.
Then the Weierstra\ss{} equation is globally minimal, the N\'eron model $E\ra\Spec(R)$ is semi-abelian,
and the closed fibers $E_a$ have Kodaira symbol $\I_v$, with $v=\val_a(\Delta)$. 
\end{lemma}

\proof
The problem is local, so we may assume that $R$ is a discrete valuation ring, say with 
uniformizer $\pi\in R$. Seeking a contradiction, we assume that the Weierstra\ss{} equation is not
minimal. Then there is a change of coordinates $x=u^2x'+r$, $y=u^3y'+usx'+t$ over the field $F=\Frac(R)$
such that the new Weierstra\ss{} equation has coefficients in $R$, with
$\val(\Delta')=\val(\Delta)-12$. Then $\val(c'_4)=\val(c_4)-4$. In particular $\Delta,c_4\in\maxid_R$, 
contradiction. Thus the equation is minimal. According to \cite{Deligne 1972}, Proposition 5.1 the closed fiber
is elliptic or multiplicative. In other words, the Kodaira symbol is $\I_v$ for some $v\geq 0$.
By the first step in the Tate Algorithm, we actually have $v=\val(\Delta)$.
\qed

\medskip
We now state two observations for Weierstra\ss{} equations and the resulting families of cubics $Y\subset\PP^2$
that are valid over arbitrary ground rings $R$.

\begin{lemma}
\mylabel{same radical ideal}
For each Weierstra\ss{} equation of the form $y^2+xy=x^3+a_2x^2+a_4x$,
the two ideals $(\Delta,c_4)$ and $(4a_2+1,a_4)$ in the ring $R$ have the same radical ideal.
\end{lemma}

\proof
It suffices to treat the case that the ground ring $R=k$ is an algebraically closed field, say  of characteristic $p\geq 0$.
According to \cite{Deligne 1972}, Proposition 5.1 the cubic $Y\subset\PP^2$ is the rational cuspidal curve if and only if the
$\Delta=c_4=0$. Using for example \cite{Magma}, we compute for our Weierstra\ss{} equation the values
\begin{equation}
\label{delta and c4}
\Delta= a_4^2 ((4a_2+1)^2-64a_4)\quadand   c_4= (4a_2+1)^2-48a_4.
\end{equation}
These two equations combined give the relation  
\begin{equation}
\label{combined relation}
\Delta = a_4^2(c_4-16a_4).
\end{equation}
Obviously, the vanishing of $4a_2+1$ and $a_4$ implies the vanishing of $\Delta$ and $c_4$.
Conversely, suppose    that $\Delta=c_4=0$. 
From the second equation in \eqref{delta and c4} we infer $p\neq 2$.
So  $a_4=0$ follows from \eqref{combined relation},  
and the second equation in \eqref{delta and c4} gives $4a_2+1=0$.
\qed

\begin{lemma}
\mylabel{unique equation}
Let $a_2,a_4,a_2',a_4'\in R$. Suppose the two Weierstra\ss{} equations
$$
y^2+xy=x^3+a_2x^2+a_4x\quadand y^2+xy=x^3+a_2'x^2+a_4'x
$$
are related by a change of coordinates that respects $P=(0,0)$.
Suppose  the ring $R$ has the property that $2\in R$ is a regular element and $R^\times \cap (1+2R)=\{\pm 1\}$.
Then the equations coincide,  and the change of coordinates is either the identity
or the sign involution $x=x'$, $y=-y'-x'$.
\end{lemma}

\proof
Write  $x=u^2x'+r$ and $y=u^3y'+su^2x'+t$. Comparing coefficients for $a_1=a_1'=1$
gives $u=1+2s$, which means $u=s=-1$ or $u=1$, $s=0$. Composing with the sign involution if necessary,
we may assume that the latter holds. Since $P=(0,0)$ is fixed, we also have $r=t=0$.
\qed

\section{From geometry to equations}
\mylabel{Geometry to equations}
 
Let $E_\QQ$ be an elliptic curve over the field of rationals $\QQ$, and  
$E\ra S$ be the N\'eron model over $S=\Spec(\ZZ)$.
Write $O\subset E$ for the zero-section,
$X\ra S$ for the minimal model, and $Y\ra S$ for the Weierstra\ss{} model.
The goal of this section is to establish the following:

\begin{theorem}
\mylabel{equation necessary}
Suppose the elliptic curve $E_\QQ$ satisfies the following three conditions:
\begin{enumerate}
\item The N\'eron model $E\ra\Spec(\ZZ)$ is semi-abelian.
\item The closed fiber $E\otimes\FF_2$ is an   elliptic curve.
\item There is an element $P\in E(\ZZ)$  of order two
whose image in  the group $E(\FF_2)$ remains non-zero.
\end{enumerate}
Then there are   integers $m,n\in\ZZ$ with $n$ odd and $\gcd(4m+1,n)=1$ such that the Weierstra\ss{} model
is given by $y^2+xy=x^3+mx^2+ nx$, and the section $P\subset Y$ is defined by the equations $x=y=0$.
Moreover, these  integers $m,n$ are unique, and we have $P\cap O=\varnothing$.
\end{theorem}

The proof is given at the end of this section, after a few preparations. Since the ring $\ZZ$ is factorial, the Weierstra\ss{} model
is a closed subscheme $Y\subset\PP^2$, and the affine scheme $Y\smallsetminus O\subset\AA^2$
is defined by a Weierstra\ss{} equation
$$
y^2+a_1xy+a_3y = x^3+a_2x^2+a_4x+a_6
$$
with coefficients $a_i\in\ZZ$. Any two such Weierstra\ss{} equations differ by a change of coordinates
$x=u^2x'+r$, $y=u^3y' + usx'+ t$, with $r,s,t\in\ZZ$ and $u=\pm1$, as described in
\cite{Tate 1972a}, Section 2.

\begin{lemma}
\mylabel{left hand side}
If the closed fiber $E\otimes\FF_2$ is an ordinary elliptic curve,
there is a Weierstra\ss{} equation with $a_1=1$ and $a_3=0$.
\end{lemma}

\proof
Let us recall from \cite{Husemoller 1987}, Chapter 3, \S6 
that up to isomorphism, there are five elliptic curves $E_1,\ldots,E_5$ over the prime field $k=\FF_2$.
The groups $E_i(k)$ are cyclic, and all   $1\leq n\leq 5$ occur as orders. Let us  choose indices in a natural fashion, so that
$|E_i(k)|=i$. To be explicit, we write
$$
E_1:\,y^2+y=x^3+x^2+1,\quad E_3:\, y^2+y=x^3\quadand E_5:\,y^2+y=x^3+x^2,
$$
which are supersingular and  have invariant $j=0$. We are more interested in
the ordinary elliptic curves, which have $j=1$ and are given by
\begin{equation}
\label{ordinary curves}
E_2:\, y^2+xy=x^3+x^2+x\quadand E_4:\, y^2+xy=x^3+x.
\end{equation}
Since we are in characteristic $p=2$ and the unit group $\FF_2^\times$ is trivial, the coefficient $a_1\in\FF_2$ of the Weierstra\ss{} equations
is invariant under every change of coordinates. For our N\'eron model $E\ra\Spec(\ZZ)$ of the elliptic curve $E_\QQ$, this means
that   $a_1 \in\ZZ$ is odd.
A change of coordinates $y=y'+sx'$ achieves $a_1=1$. A further change of coordinate $y=y'+t$
then yields $a_3=0$ as well. 
\qed

\begin{proposition}
\mylabel{also a6}
If conditions (ii) and (iii) from Theorem \ref{equation necessary} hold, one may choose a Weierstra\ss{} equation 
so that $P\subset Y$ is given by $P=(0,0)$, and the coefficients satisfy $a_1=1$ and $a_3=a_6=0$. 
Furthermore,   $a_4\in\ZZ$ then must be odd.
\end{proposition}

\proof
We saw in the previous proof that the supersingular elliptic curves over $\FF_2$ contain no element of order two,
hence we may find a Weierstra\ss{} of the form $y^2+xy=x^3+a_2x^2+a_4x+a_6$. This form is preserved  under each  change of coordinates
$x=x'-2t$, $y=y'+t$. Since the subscheme $P\subset Y$ is disjoint from the zero-section, it has coordinates $P=(m,n)$
for some integers $m,n\in\ZZ$. So after a change of coordinates, we may assume $n=0$.
The sign involution on $E$ is given by $-(x,y)=(x,-y-x)$. For the element $P\in E(\ZZ)$ of order two, this means
$(m,0) = (m,-m)$, thus $m=0$.
Summing up, $x=y=0$ is a solution for our Weierstra\ss{} equation, and this means $a_6=0$.
Seeking a contradiction, we suppose that $a_4$ is even. Then the fiber $Y\otimes\FF_2$ is given by
$y^2+xy=x^3+x^2$ or $y^2+xy=x^3$. Both curves are singular at the point $(0,0)$, contradiction.
\qed

\begin{proposition}
\mylabel{gcd condition}
Suppose   $Y$ is given by a Weierstra\ss{} equation of the form $y^2+xy=x^3+mx^2+nx$. Let $p>0$ be a  prime.
Then the fiber $E\otimes \FF_p$ is additive if and only if we have $p\mid \gcd(4m+1,n)$.
\end{proposition}

\proof
Consider the ideal $\ideala=(4m+1,n)$.
According to Lemma \ref{same radical ideal}, the fiber  $E\otimes\FF_p$ is additive if and only if $\ideala\subset (p)$,
or equivalently $p\mid\gcd(4m+1,n)$.
\qed

\medskip
\emph{Proof of Theorem \ref{equation necessary}:}
Suppose we have an elliptic curve $E_\QQ$ whose N\'eron model $E\ra\Spec(\ZZ)$ satisfies conditions (i)--(iii).
According to Proposition \ref{also a6}, we may choose a Weierstra\ss{} equation of the form $y^2+xy=x^3+mx^2+nx$.
In light of Proposition \ref{gcd condition} we have $\gcd(4m+1,n)=1$.
The coefficients $m,n\in\ZZ$ are unique by Lemma \ref{unique equation}.
The Weierstra\ss{} model $Y\subset\PP^2$ is the zero-locus for the homogeneous equation
$
Y^2Z+XYZ = X^3+mX^2Z+nXZ^2
$.
The images of $O$ and $P$ on the Weierstra\ss{} model
are given by the equations   $X=Z=0$ and $X=Y=0$, which are clearly disjoint. 
Thus $O,P\subset E$ are disjoint.
\qed

\section{Analysis of the Weierstra\ss{} equation}
\mylabel{Analysis}

Let $m,n\in\ZZ$ be two integers with  $n$   odd and $\gcd(4m+1,n)=1$.
We now consider the family of cubics $Y\subset\PP^2$ defined by the Weierstra\ss{} equation
\begin{equation}
\label{weierstrass equation}
y^2+xy=x^3+mx^2+nx.
\end{equation}
The discriminant   $\Delta=n^2 ((4m+1)^2-64n)$ is odd, in particular the      generic fiber
$E_\QQ=Y_\QQ$ is an elliptic curve. Combining Lemma \ref{same radical ideal} and Proposition  \ref{globally minimal}, we obtain:

\begin{proposition}
\mylabel{globally minimal weierstrass model}
The   Weierstra\ss{} equation \eqref{weierstrass equation} is globally minimal, and
the family of cubics $Y\subset\PP^2$ coincides with the Weierstra\ss{} model for the elliptic curve $E_\QQ$.
\end{proposition}
 
Let $E\ra\Spec(\ZZ)$ be the N\'eron model, with zero-section $O\subset E$.
The section of the Weierstra\ss{} model $Y$ given by $x=y=0$ induces 
a section $P\subset E$. The following is an immediate converse for Theorem \ref{equation necessary}:

\begin{proposition}
\mylabel{equation sufficient}
The relative group scheme $E\ra\Spec(\ZZ)$ is semi-abelian,
the element $P\in E(\ZZ)$ has order two, and the section $P\subset E$ is disjoint from the zero-section.
The closed fiber $E\otimes\FF_2$ is an ordinary elliptic curve, which contains two rational points if $m$ is odd,
and four rational points if $m$ is even.
\end{proposition}

\proof
The sign involution on $E$ is given by $(x,y)\mapsto (x,-y-x)$, hence the element $P\in E(\ZZ)$ has order two.
Obviously, the images $(0:1:0)$ and $(0:0:1)$  on the Weierstra\ss{} model $Y\subset\PP^2$  of the sections
$O,P\subset E$ are disjoint,
and the closed fiber $E\otimes\FF_p$ belongs to the ordinary elliptic curves in \eqref{ordinary curves}.
The N\'eron model $E$ is semi-abelian by Proposition \ref{gcd condition}. Obviously, the close fiber $E\otimes\FF_2$
occurs in our list of ordinary elliptic curves \eqref{ordinary curves}, and the assertion on $E(\FF_2)$ is immediate.
\qed

%
 
\begin{proposition}
\mylabel{component group}
Let $p>0$ be a prime. The element $P\in E(\ZZ)$ has non-trivial class in the component group scheme $\Phi_p$ 
if and only if $p\mid n$.
\end{proposition}

\proof
The image of $P\in E(\ZZ)$ in the group   $\Phi_p(\FF_p)$ vanishes if and only if the image 
of the  subscheme $P\subset E$ 
on the Weierstra\ss{} model $Y$ is disjoint from the singular locus $\Sing(Y\otimes\FF_p)$.
The structure morphism $Y\ra\Spec(\ZZ)$ is smooth on some common open neighborhood of  the zero-section $O\subset Y$ and 
the elliptic curve $Y\otimes\FF_2$.
Hence the locus of non-smoothness $\Sing(Y/\ZZ)$ is the closed subscheme in $\AA^2\otimes\ZZ[1/2]$
defined by the Weierstra\ss{} equation \eqref{weierstrass equation} and its partial derivatives
\begin{equation}
\label{partial derivatives}
2y+x=0\quadand y=3x^2+2mx+n.
\end{equation}
For later use, we record that $\Sing(Y/\ZZ)$ therefore equals the spectrum of the ring
\begin{equation}
\label{locus non-smoothness}
\ZZ[1/2,x,y]/( 2y+x, -8y^3+(4m+1)y^2-2ny, 12y^2-(4m+1)y+n),
\end{equation}
which can be seen by substituting $x=-2y$ in the other two equations.
By abuse of notation, we now identify the section $P\subset E$ with its
image in the Weierstra\ss{} model $Y$, where it is is given by the equation $x=y=0$.
Thus the intersection $P\cap \Sing(Y/\ZZ)$ is defined by the equations \eqref{weierstrass equation} and \eqref{partial derivatives},
together with $x=y=0$. Clearly, this is the spectrum of $\ZZ[1/2]/(n)$, which coincides with
$\ZZ/n\ZZ$  because $n$ is odd.
 This shows that $P$ intersects $\Sing(Y\otimes\FF_p)$ if and only if $p\mid n$. The assertion follows.
\qed

\medskip
As an immediate consequence, we get:

\begin{corollary}
\mylabel{narrow}
The element $P\in E(\ZZ)$ is narrow if and only if $n=\pm 1$.
\end{corollary}
 
\medskip
Fix some prime $p>0$. The fiber $E\otimes\FF_p$  has Kodaira symbol $\I_v$ for some integer $v\geq 0$.
Let us make some general observations on the  multiplicative case
$v\geq 1$.   The    component of the origin $G^0$ of the algebraic group scheme 
$G=E\otimes\FF_p$ is either the multiplicative group
$\GG_{m,\FF_p}$ or its quadratic twist $\widetilde{\GG}_{m,\FF_p}$ with respect the involution $\lambda\mapsto 1/\lambda$ 
and  the cyclic extension $\FF_p\subset\FF_{p^2}$.  
For the fibers of the  Weierstra\ss{} model, this means the following:

\begin{proposition}
Suppose    $v\geq 1$ and   $G^0=\widetilde{\GG}_{m,\FF_p}$. Then the fiber $Y\otimes\FF_p$ is obtained from the projective line $\PP^1\otimes\FF_p$  
by replacing an $\FF_{p^2}$-valued point by some rational point.
\end{proposition}

\proof
Set $C=Y\otimes\FF_p$ and let $\nu:C'\ra C$ be the normalization map. Then $C'$ is the projective line over $\FF_p$, because this field is perfect and
has trivial Brauer group. Let $B\subset C$ be the support of the coherent sheaf $\nu_*(\O_{C}')/\O_C$ and $B'\subset C'$ be its preimage.
Then  there is an exact sequence $0\ra \O_C\ra\O_{C'}\oplus\O_B\ra\O_{B'}\ra 0$, compare \cite{Schroeer; Siebert 2006}, Section 4.
The long exact sequence reveals that $h^0(\O_B)=1$ and $h^0(\O_{B'})=2$. Since the fiber is multiplicative, the scheme $B'$ is \'etale.
So it is either the spectrum of $A=\FF_p\times\FF_p$ or $A=\FF_{p^2}$.
In turn, the number of rational points in $G^0=C\smallsetminus B=C'\smallsetminus B'$ is either $n=(p+1)-2$ or $n=(p+1)$.
The assertion follows from the fact that  $\GG_m(\FF_p)=\FF_p^\times$ has order $p-1$.
\qed

\medskip
For the regular model $X$, the situation is as follows:
 
\begin{proposition}
Suppose    $v\geq 2$ and $G^0=\widetilde{\GG}_{m,\FF_p}$.
Then $C=X\otimes\FF_p$ is  a chain of curves $C_0,C_1,\ldots,C_w$, where the component $C_0$ of the origin is a projective
line over $\FF_p$, and the    $C_1,\ldots,C_{w-1}$ are projective lines over $\FF_{p^2}$.
The   intersections are 
$C_i\cap C_{i+1}\simeq\Spec(\FF_{p^2})$ for $0\leq i\leq w-1$. 
If $v$ is odd, we have $w=(v+1)/2$ and $C_w=\PP^1\otimes \FF_p$.
If $v$ is even, we have  $w=(v+2)/2$  and  $C_w$ arises from $\PP^1\otimes\FF_{p^2}$ by replacing an $\FF_{p^2}$-valued point
by some rational point.
\end{proposition}

\proof
According to the Tate Algorithm, the base-change $D=X\otimes\FF_{p^2}$ is a cycle   of   smooth rational curves over $\FF_{p^2}$, say with
irreducible components $D_0,\ldots,D_{v-1}$ and intersections $D_j\cap D_{j+1}=\Spec(\FF_{p^2})$,
where we regard   indices as congruence classes modulo $v$.
The base-change comes with an  action of the Galois group $\Gamma=\Gal(\FF_{p^2}/\FF_p)$, which is cyclic of order two,
and the quotient is $C=X\otimes\FF_p$.  
By assumption, the  generator $\sigma\in\Gamma$ stabilizes
$D_0=C_0\otimes\FF_{p^2}=\Spec\FF_{p^2}[t] \cup \Spec\FF_{p^2}[t^{-1}]$, with induced action given by $t\mapsto t^{-1}$.
It follows inductively   that $\sigma(D_i)=D_{-i}$ for $0\leq i\leq v-1$.
 
If $v\geq 2$ is even, the $\Gamma$-action on the irreducible components of $D$ reveals that  $w=1+(v-2)/2+1=(v+2)/2$. 
We have $D_i\cap D_{-i}=\varnothing$  for $1\leq i\leq w-1$,
hence the projections $D_i\ra C_i$ are isomorphisms. The $\Gamma$-action on $D_w$ is free, so $C_w=D_w/\Gamma=\PP^1\otimes\FF_p$.
If $v\geq 3$ is odd, then $w=1+(v-1)/2=(v+1)/2$. Now the projections $D_i\ra C_i$ are isomorphisms for $1\leq i\leq w-1$.
In contrast, $D_w\ra C_w$ is only birational, and identifies the $\FF_{p^2}$-valued point $D_w\cap D_{w+1}$ to a rational point.
Finally, one easily computes $C_i\cap C_{i+1}\simeq\Spec(\FF_{p^2})$.
\qed


\medskip
If $G^0=\GG_{m,\FF_p}$, one says that $E_\QQ$ has \emph{split multiplicative reduction},
and   the corresponding fibers of $E$, $X$ and $Y$ are called \emph{untwisted}.
On the other hand, if $G^0=\widetilde{\GG}_{m,\FF_p}$, one says that the elliptic curve
has \emph{non-split multiplicative reduction}, and the corresponding fibers are \emph{twisted}.

Let us now unravel what this concretely means for  our Weierstra\ss{} equation \eqref{weierstrass equation}.
Here the Legendre symbol $(\frac{d}{p}) =\pm 1$ is useful, which gives the class of of an integer $d$ in
$\FF_p^\times/\FF_p^{\times 2}$, where $p>0$ is an odd prime not dividing $d$.

\begin{proposition}
\mylabel{fiber}
For each prime $p>0$,  the fiber $E\otimes\FF_p$ has  Kodaira symbol $\I_v$ with index
$v=2\val_p(n)+\val_p ((4m+1)^2-64n))$.
Moreover, in the multiplicative case $v\geq 1$, $p\neq 2$ the following holds:
\begin{enumerate}
\item
If $p\mid n$, the singularity on $Y\otimes \FF_p$ has coordinates  $x=y=0$, and the
fiber $E\otimes\FF_p$ is untwisted if and only if   $(\frac{4m+1}{p})=1$.
\item
If $p\nmid n$, the singularity on $Y\otimes\FF_p$ has coordinates $x=-(4m+1)/8$ and $y=(4m+1)/16$, and the fiber
is untwisted if and only if the Legendre symbol satisfies   $(\frac{4m+1}{p})=(-1)^{(p^2+4p-5)/8}$.
\end{enumerate}
\end{proposition}

\proof
Our Weierstra\ss{} equation has discriminant $\Delta=n^2 ((4m+1)^2-64n)$.  
The   assertion on the index $v\geq 0$ follows from the minimality of the Weierstra\ss{} equation and the Tate Algorithm.
Now suppose that $v\geq 1$, such that $p\neq 2$. The change of coordinates $x=x'$ and $y=y'-\frac{1}{2}x'$ yields a
Weierstra\ss{} equation in simplified form
$y^2=x^3+ \frac{4m+1}{4} x^2 +n x$.
In case (i), the closed fiber $Y\otimes\FF_p$ has an affine part given by the ring 
$$
A=\FF_p[x,y]/(y^2- x^2(x+\frac{4m+1}{4})).
$$
Clearly, the singularity has coordinates $x=y=0$.
The fraction $t=2y/x$ lies in the normalization $A\subset A'$, and satisfies the relation $t^2=4x+(4m+1)$.
The fiber ring over the singularity is $L=k[t]/(t^2-(4m+1))$, and we have $E\otimes\FF_p=\PP^1\otimes\FF_p\smallsetminus\Spec(L)$.
The ring $L$ a field if and only if 
$4m+1\in\FF_p^\times$ is not a square, and (i) follows.

Now suppose we are in case (ii), such that  $p$ divides $(4m+1)^2-64n$. Over $\FF_p$, we get
$x^2+ \frac{4m+1}{4} x +n = (x+\frac{4m+1}{8})^2$. In turn, the singular locus of $Y\otimes\FF_p$
has coordinates $x'=-\frac{4m+1}{8}$, $y'=0$ in the new coordinates, and $x=-\frac{4m+1}{8}$, $y=\frac{4m+1}{16}$
in the original coordinates. Arguing as in the previous paragraph, the element $t=y/(x+\frac{4m+1}{8})$ lies
in the normalization, and the fiber ring becomes $L=k[t]/(t^2+\frac{4m+1}{8})$.
This is a field if and only if $-2(4m+1)\in \FF_p^\times$ is not a square.
According to \cite{Serre 1973}, Chapter I, \S3, Theorem 5, we have $(\frac{-2}{p}) = (-1)^d$ with exponent $d=(p-1)/2 + (p^2-1)/8 = (p^2+4p-5)/8$,
and (ii) follows.
\qed

\medskip
Let us say that the fiber  $E\otimes\FF_p$ has \emph{even Kodaira symbol} if the Kodaira symbol
is $\I_v$ for some even integer $v\geq 0$. This means that either $E\otimes\FF_p$ is an elliptic curve,
or that the component group scheme $\Phi_p$ has even order.  
For each integer $r\geq 1$, write  $E[r]$ for the kernel of the multiplication map $r:E\ra E$, which is defined by
the cartesian diagram
$$
\begin{CD}
E[r]	@>>>	\Spec(\ZZ)\\
@VVV		@VVO V\\
E	@>>r>	E.
\end{CD}
$$
Note that the formation   of    kernels   commutes with base-change, because fiber products do so. Moreover,
$E[r]$ inherits the structure of relative group scheme, whose structure morphism is separated and of finite type.

We now come to our first main result. Recall that $m$ and $n$ are integers, with $n$ odd and relatively prime to $4m+1$.
 
\begin{theorem}
\mylabel{additional 2-division}
For the N\'eron model $E\ra\Spec(\ZZ)$ of  $y^2+xy=x^3+mx^2+nx$,
and the section $P$ defined by $x=y=0$,
the following are equivalent:
\begin{enumerate}
\item There is another element $Q\neq P$ of order two in $E(\ZZ)$.
\item The relative group scheme $E[2]$ is isomorphic to $\ZZ/2\ZZ\oplus\mu_2$.
\item For every prime $p>0$, the fiber $E\otimes\FF_p$ has   even  Kodaira symbol.
\item We have $n = -d(4m +1+16d)$ for some odd   $d$ with $\gcd( 4m+1,d)=1$.
\end{enumerate}
In this situation, the three elements of $E(\ZZ)$ of order two are given by
$$
P=(0,0),\quad Q=(4d,-2d)\quadand P+Q=(-\frac{4m+1+16d}{4}, \frac{4m+1+16d}{8}).
$$
\end{theorem}
 
\proof
The implication (ii)$\Rightarrow$(i) is trivial. For (iv)$\Rightarrow$(iii) we compute
$$
(4m+1)^2-64n= (4m+1)^2+64d(4m +1+16d) =(4m + 1+ 32d)^2.
$$
In light of \eqref{delta and c4} the discriminant of our Weierstra\ss{} equation becomes
\begin{equation}
\label{square decomposition}
\Delta= d^2(4m +1+16d)^2(4m + 1+ 32d)^2,
\end{equation}
and we conclude with Proposition \ref{fiber} that all closed fibers $E\otimes\FF_p$ have even
Kodaira symbol.  

We now show (i)$\Rightarrow$(iv).
Let $Q\neq P$ be another element of order two.
Replacing $Q$ by $Q+P$ if necessary, we may assume that its image  in the group $E(\FF_2)$
is non-zero. According to Theorem \ref{equation necessary}, the family of cubics $Y\subset\PP^2$ admits
a description with another Weierstra\ss{} equation
$$
y'^2+x'y'=x'^3+m'x'^2+n'x'
$$
for some unique integers $m',n'\in \ZZ$ with $n'$ odd and $\gcd(4m'+1,n')=1$, where in the new coordinates
the closed subscheme $Q\subset Y$ is defined by $x'=y'=0$.
The two Weierstra\ss{} equations are related by a change of variables $x=u^2x'+r$,
$y=u^3y'+su^2x'+t$ with $u=\pm 1$ and $r,s,t\in\ZZ$. Composing with the sign involution if necessary, we may assume that $u=1$.
We now use  \cite{Tate 1972a}, equations (2.2) to compare coefficients.
The condition $ a_1'=a_1=1$ yields $s=0$, and   $a'_3=a_3=0$ gives
$r=-2t$. In turn, we get 
$$
n'=a_4'=a_4+2ra_2-t+3r^2=n-4tm-t+ 12t^2. 
$$
In particular, $t$ is a non-zero even integer, which we write as $t=-2d$. 
Comparing coefficients at $a_6'=a_6=0$ gives
$$
0=ra_4 + r^2a_2 + r^3 - t^2 - rt = 4dn + 16d^2m +64d^3  + 4d^2.
$$
Dividing by $4d\neq 0$ yields the desired  equation $n+4dm +16d^2+ d=0$.
The integer $d$ must be odd, because $n$ is odd.

We now come to the implication (iii)$\Rightarrow$(ii). 
To simplify notation, write $G=E[2]$ and $S=\Spec(\ZZ)$.
Since $E$ is semi-abelian,
the multiplication map $r:E\ra E$ is quasi-finite and flat for all integers $r\neq 0$, according to 
\cite{Bosch; Luetkebohmert; Raynaud 1990}, Section 7.2, Lemma 2.
Being a  base-change, the structure morphism    $G\ra S$ is quasi-finite and flat.
Each geometric fiber for  $E\ra S$ is either an elliptic curve, or   $\GG_m\times\ZZ/v\ZZ$
with an even integer $v\geq 1$.  It follows that the fibers $G\otimes\FF_p$ are finite group schemes of constant order four.
By Zariski's Main Theorem,   $G$ is an open subscheme of some finite $\ZZ$-scheme $W$.
Replacing $W$ by the closure of $G$ me may assume that $G\subset W$ is schematically dense.
Then the structure morphism $W\ra S$, which  is finite and flat, has $\deg(W/S)=4$.
Hence for each prime $p>0$, the inclusion $G\otimes\FF_p\subset W\otimes\FF_p$ is an equality. The Nakayama Lemma gives
$G=W$, so the structure morphism $G\ra\Spec(\ZZ)$ is finite and flat of degree four.
 
The section $P\subset E$ yields a homomorphism $\ZZ/2\ZZ\ra E[2]$ of group schemes over $S$. 
This is a monomorphism, because the intersection $P\cap O$ is empty.
The quotient $H$ exists as a  group scheme over $\ZZ$, its formation commutes with base-change,
and the structure morphism $H\ra S$ is free of rank two.
By the Tate--Oort Classification (\cite{Tate; Oort 1970}, Corollary to Theorem 3), we either have  $H=\mu_2$ or $H=\ZZ/2\ZZ$.
In the latter case, the group $G\otimes\FF_2$ would be \'etale, contradicting the fact that multiplication by $r=2$
on the elliptic curve $E\otimes\FF_2$ is not \'etale.
Thus we have a short exact sequence $0\ra \ZZ/2\ZZ\ra G\ra\mu_2\ra 0$. By Lemma \ref{minkowski} below,
such extensions split, such that  (ii) follows.


It remains to verify the coordinates for the   two-torsion sections,
assuming that the equivalent conditions (i)--(iv) hold.
By definition we have $P=(0,0)$.
In light of  the equation $n = -d(4m+1+16d)$ we make the
change of coordinates
$$
x=x'+4d\quadand y=y' -2d,
$$
which transforms the old Weierstra\ss{} equation into $y'^2+x'y'=x'^3+m'x'^2+n'x'$, where
\begin{equation}
\label{new coefficients}
m'=m+12d\quadand n'=  n+8dm+2d+ 48d^2 = d(4m+1+32d).
\end{equation}
The two-torsion section $x'=y'=0$ in the new Weierstra\ss{} equation corresponds to the
two-torsion section $Q\subset Y$ defined by $x=4d$ and $y=-2d$.
In other words, $Q=(4d,-2d)$. The coordinates for $P+Q$ follow from the group law for   elliptic curves.
\qed

\medskip
In the preceding proof, we have used a special case of the following observation:

\begin{lemma}
\mylabel{minkowski}
Let $N$ be a finite group, and write $N_S$ for the corresponding constant group scheme  over $S=\Spec(\ZZ)$.
Then every extension  $1\ra N_S\ra G\ra \mu_2\ra 0$ of group schemes splits, and the resulting semidirect product $G=N_S\rtimes \mu_2$
is actually a direct product.
\end{lemma}

\proof
The affine group scheme $\mu_2$ is the spectrum of the ring $\ZZ[T]/(T^2-1)$, which sits in a cartesian diagram
$$
\begin{CD}
\FF_2\times\FF_2	@<<<	\ZZ\times\ZZ\\
@AAA				@AAA\\
\FF_2			@<<<	\ZZ[T]/(T^2-1).
\end{CD}
$$
Roughly speaking, the scheme $\mu_2 $ is a connected union $S_1\cup S_2$ of two copies of $S$, where the  two copies of $\Spec(\FF_2)$ are identified.

According to Minkowski's Theorem  (\cite{Neukirch 1999}, Chapter III, Theorem 2.17), every number field $\QQ\subset K$ has discriminant $d_K\neq \pm1$.
It follows that 
the scheme $S=\Spec(\ZZ)$ is simply-connected,
in the sense that for each  finite \'etale morphism $S'\ra S$, the scheme $S'$ is   the  disjoint union of finitely many copies of $S$.
In other words, the algebraic fundamental group $\pi_1^\alg(S,a)$ vanishes (\cite{SGA 1}, Expos\'e V),
for any geometric point $a:\Spec(\Omega)\ra S$.

The projection $G\ra\mu_2$ is finite and \'etale, being a torsor for the finite \'etale group   scheme $N_{\mu_2}$ over $\mu_2$.
Hence the same holds for the base-changes $G\times_{\mu_2}S_i\ra S_i$. The structure maps $S_i\ra S$ are isomorphisms, so the $S_i$ are simply-connected,
and  $G\times_{\mu_2}S_i$ are just disjoint unions of copies of $S_i$. 
Note that this property may fail  for  other number rings $R\neq \ZZ$.

The neutral element $e\in G\otimes\FF_2$ is contained in $G\times_{\mu_2}S_i$.
Let  $s_i:S_i\ra G\times_{\mu_2}S_i$ be the unique section corresponding to the copy of $S_i$ passing through this point.
Then $s_1$ and $s_2$ coincide on $S_1\cap S_2=\Spec(\FF_2)$, and thus define
a section $s:\mu_2\ra G$.
The expression $\psi(x,y)=s(xy)/s(x)s(y)$, where we write the group laws   in multiplicative fashion,
defines a morphism of schemes $\psi:\mu_2\times\mu_2\ra N_S$.
It factors over   $\left\{n\right\}\times S\subset N_S$ for some group element $n\in N$, because both   $\mu_2\times\mu_2$ and $S$ are
connected. Base-changing to $\FF_2$, we see that $n\in N$ must be the neutral element.
In turn, the section of schemes $s:\mu_2\ra G$  is a homomorphism of group schemes, and thus 
splits  the extension of group schemes.

Consequently, the extension is given by a semidirect product $G=N_S\rtimes \mu_2$, formed with respect to some
homomorphism $\varphi:\mu_2\ra\Aut(N_S)=\Aut(N)_S$. Arguing as above, we see that $\varphi$ is trivial,
hence $G=N_S\times\mu_2$. 
\qed

\section{Quotients by 2-division points}
\mylabel{Quotients}

We keep the notation as in the previous section, and study the N\'eron model $E$
for the Weierstra\ss{} equation $y^2+xy=x^3+mx^2+nx$, where $n=-d(4m+1+16d)$ for some odd integer $d$ with $\gcd(4m+1,d)=1$.
We saw in Theorem \ref{additional 2-division} that there are two elements $P,Q\in E(\ZZ)$
of order two whose images in $E(\FF_2)$ keep order two.
Moreover, the sum $P+Q$ generates a copy of the multiplicative group scheme $\mu_2\subset E$.
The quotient $E/\mu_2$ is a relative group scheme over the ring $\ZZ$ whose structure morphism
is separated and of finite type. It comes with a canonical section of order two that is disjoint from the zero-section,
because $E[2]/\mu_2=\ZZ/2\ZZ$. From this one may deduce that Theorem \ref{equation necessary} applies
to the generic fiber $(E/\mu_2)_\QQ$. In fact, it is possible to  infer directly the  equation
for its Weierstra\ss{} model:

\begin{proposition}
\mylabel{quotient model}
The   elliptic curve $(E/\mu_2)_\QQ$ is given by the globally minimal Weierstra\ss{} equation
$y^2+xy=x^3+ (m+6d)x^2+d^2x$.
\end{proposition}
 
\proof
Applying Velu's Formula \cite{Velu 1971} with the subgroup in $(E/\mu_2)_\QQ$ generated by the  2-division point
$$
P+Q=(-\frac{4m+1+16d}{4}, \frac{4m+1+16d}{8}),
$$
we get a Weierstra\ss{} equation $y^2+xy=x^3+mx^2+a_4x^2+a_6$ with new coefficients
\begin{gather*}
a_4=-5m^2 + 64md - \frac{5}{2}m - 176d^2 + 16d - \frac{5}{16},\\
a_6=3m^3 - 64m^2d + \frac{9}{4}m^2 + 432md^2 - 32md + \frac{9}{16}m - 896d^3 +     108d^2 - 4d + \frac{3}{64}.
\end{gather*}
A change of coordinates $x=x'-2t$, $y=y'+t$ with $t= -\frac{1}{2}m - 4d - \frac{1}{8}$
transforms this into
$y^2 + xy = x^3 + (4m + 24d + \frac{3}{4})x^2 + 16d^2x$.
A further change of coordinates with $x=x'$, $y=y'+\frac{1}{2}x'$ gives $y^2+2xy=x^3+(4m + 24d)x^2 + 16d^2x$,
which immediately leads to the desired Weierstra\ss{} equation.
\qed

\medskip
The elements $P,Q\in E(\ZZ)$ may pass through different components of the fibers $E\otimes\FF_p$.
We may easily unravel the situation, by considering their difference, which is also their sum:

\begin{lemma}
\mylabel{non-trivial in component group}
The element $P+Q\in E(\ZZ)$ 
has non-trivial class  in the component group scheme $\Phi_p$ if and only if the prime $p$ divides $(4m+1+16d)(4m+1+32d)$.
\end{lemma}

\proof
The assertion is trivial for $p=2$. Assume now that $p$ is odd.
Clearly, the element $P+Q$ has non-trivial class in $\Phi_p$ if and only if its image on the Weierstra\ss{} model passes through the
singularity of the fiber $Y\otimes\FF_p$. This image has coordinates $x=-(4m+1+16d)/4$ and $y=(4m+1+16d)/8$   
by Theorem \ref{additional 2-division},
and the singular locus is given by $12y^2-(4m+1)y+n=0$ and $x=-2y$, according to   \eqref{locus non-smoothness}.
We see $P+Q\not\equiv 0$ in $\Phi_p$ if and only if $p$ divides 
$$
16m^2 + 192md + 8m + 512d^2 + 48d + 1 = (4m+1+16d)(4m+1+32d),
$$
and the result follows.
\qed

\medskip
The translation action of $\mu_2$ on the N\'eron model $E$ extends to an action of the minimal model $X$.
The induced action on the fibers  $X\otimes\FF_p$ may or may not be free.
 
\begin{proposition}
\mylabel{action non-free}
The induced action of $\mu_2\otimes\FF_2$ on the fiber $X\otimes\FF_p$ is not free if and only if $p\mid d$.
In this case, the locus of fixed points coincides with $\Sing(X\otimes\FF_p)$. Moreover,
the minimal resolution of singularities for the quotient scheme $X/\mu_2$ is the minimal model
of the elliptic curve $(E/\mu_2)_\QQ$.
\end{proposition}

\proof
The action is non-free if and only if $E\otimes\FF_p$ is multiplicative, and   $P+Q\equiv O$ in the component group scheme $\Phi_p$.
In this case, the $\mu_2$-action stabilizes every irreducible component and fixes precisely the singularities on $X\otimes\FF_p$.
According to equations  \eqref{square decomposition}, 
the discriminant $\Delta$ of our  Weierstra\ss{} equation is the square of 
$ d(4m+1+16d)(4m+1+32d)$.
Thus the fiber $E\otimes\FF_p$ is multiplicative if and only if $p$ divides this number.
By Lemma \ref{non-trivial in component group}, 
we have $P+Q\equiv 0$ in the component group  scheme $\Phi_p$ if and only if $p$ does not divide $(4m+1+16d)(4m+1+32d)$.
The first assertion follows.
 
The fixed points $a\in X$ for the $\mu_2$-action yield rational double points of type $A_1$ on the quotient scheme $X/\mu_2$.
Let $X'\ra X/\mu_2$ be the resolution of singularities. The three schemes in $X\ra X/\mu_2\leftarrow X'$ are Gorenstein.
The dualizing sheaf $\omega_{X/\ZZ}$ is trivial, hence the dualizing sheaves on $X/\mu_2$ and $X'$ are numerically trivial.
It follows that $X'\ra\Spec(\ZZ)$ contains no $(-1)$-curve, thus $X'$ is the minimal model of $(E/\mu_2)_\QQ$.
\qed

\section{Classification of curves with additional sections}
\mylabel{Classification}

We now classify elliptic curves $E_\QQ$ whose N\'eron model $E$ is semi-abelian
and whose Mordell--Weil group $E(\ZZ)$ has certain elements of order two.
Our first   result is:

\begin{theorem}
\mylabel{classification order two}
Up to isomorphism, there  are exactly two elliptic curves $E_\QQ$ such that its N\'eron model $E$ has
the following properties:
\begin{enumerate}
\item The structure morphism $E\ra\Spec(\ZZ)$ is semi-abelian.
\item The closed fiber $E\otimes\FF_2$ is an   elliptic curve.
\item There is a narrow element $P\in E(\ZZ)$ of order two whose images in
$E(\FF_2)$ are non-zero, and another element $Q\neq P$ of order two.
\end{enumerate}
These elliptic curve are given by the global minimal Weierstra\ss{} equation
$$
y^2+xy = x^3-4x^2-x\quadand y^2+xy=x^3+4x^2+x,
$$ 
The former has invariant $j=20346417/289$ and the only singular fiber occurs at $p=17$. 
The latter has $j=13997521/225$, with singular fibers at $p=3$ and $p=5$. All these singular fibers have Kodaira symbol   $\I_2$. 
\end{theorem}

\proof
Suppose $E_\QQ$ satisfies the conditions (i)--(iii).
According to Theorem \ref{equation necessary}, the Weierstra\ss{} model $Y$ is given by the equation
$y^2+xy=x^3+mx^2+nx$ where $n$ is odd and $\gcd(4m+1,n)=1$, and the element $P\in E(\ZZ)$ is given
by $x=y=0$. Since this element is narrow, we must have $n=\pm 1$, by Corollary \ref{narrow}.
Since there is another 2-division point $Q\neq P$, Theorem \ref{additional 2-division} tells us that
$n=-d(4m+1+16d)$ for some integer $d$. We thus have $d=\pm 1$ and $4m+16d=0$, thus   $m=-4d$.
Consequently the only solutions are $n=-1$, $d=1$, $m=-4$ and $n=1$, $d=-1$, $m=4$.
This gives our two Weierstra\ss{} equations.

Conversely, we have to check that the two equations have the stated properties.
The equations are relatively minimal by Proposition
\ref{globally minimal weierstrass model}, and the N\'eron model $E$ satisfies conditions (i)--(iii) 
according to Proposition \ref{equation sufficient}.
The description of the singular fibers follow from Proposition \ref{fiber}.
One may compute their $j$-invariant  with \cite{Magma}.
\qed

\medskip
Let $E$ be the N\'eron model for the Weierstra\ss{} equation $y^2+xy=x^3 + 4nx^2+nx$ with $n=\pm 1$
as in the preceding theorem.
After replacing $Q$ by $Q+P$, we may assume that both $P,Q\in E(\ZZ)$ stay non-trivial in $E(\FF_2)$.
Consider the subgroup scheme $\mu_2\subset E$ generated by the sum $P+Q$. 
The induced $\mu_2$-action on the minimal model $X$ is free, according to Proposition \ref{action non-free},
and   the quotient $X/\mu_2$ arises from the Weierstra\ss{} equation $y^2+xy=x^3\pm 2x^2+x$, by Proposition \ref{quotient model}.
It is somewhat surprising that these two curves can be characterize in terms of 4-division points:

\begin{theorem}
\mylabel{classification order four}
Up to isomorphism, there  are exactly two elliptic curves $E_\QQ$ such that its N\'eron model $E$ has
the following properties:
\begin{enumerate}
\item The structure morphism $E\ra\Spec(\ZZ)$ is semi-abelian.
\item The fiber $E\otimes\FF_2$ is an ordinary elliptic curve. 
\item There is a  narrow element $R\in E(\ZZ)$ of order four    whose image in $E(\FF_2)$ keeps order four. 
\end{enumerate} 
These elliptic curve are given by the global minimal Weierstra\ss{} equation
$$
y^2+xy=x^3+2x^2+x\quadand y^2+xy = x^3-2x^2+x.
$$ 
The former has invariant $j=35937/17$, and the only singular fiber appears at $p=17$.
The latter has $j=-1/15$, with singular fibers at $p=3$ and $p=5$. 
All these singular fibers have Kodaira symbol   $\I_1$
\end{theorem}

\proof
Suppose first that $E_\QQ$ satisfies conditions (i)--(iii).
The element $P=2R$  in $E(\ZZ)$ is narrow, has order two,  and remains non-zero in $E(\FF_2)$.
Consider the group $G=\ZZ/2\ZZ$. It acts via translation by $P$
on the N\'eron model $E$, with induced actions on the minimal model $X$ and the Weierstra\ss{} model $Y$.
The $G$-action on is free on the open subset $E\subset X$, and fixes the complement $\Sing(X/\ZZ)=X\smallsetminus E$.
This complement is finite, because $E\ra\Spec(\ZZ)$ is semi-abelian, and disjoint from the fiber $X\otimes\FF_2$,
because $E\otimes\FF_2$ is an elliptic curve. 
It follows that for each point $a\in\Sing(X/\ZZ)$,
the image $b\in X/G$ yields a rational double point of type $A_1$. Let $X'\ra X/G$ be the minimal resolution of singularities.
As in the proof for Proposition \ref{action non-free}, one sees that $X'$ is the minimal model of the elliptic curve $(E/G)_\QQ$.
Write $E'$ for its  N\'eron model.

We see that $E'\ra\Spec(\ZZ)$ is semi-abelian, and 
for each closed fiber $E\otimes\FF_p$ with Kodaira symbol $\I_v$, the closed fiber $E'\otimes\FF_p$ has Kodaira symbol
$\I_{2v}$. Moreover, the image $P'\subset E'$ of $R\subset E$, which is the quotient of the closed subscheme 
$R\cup (R+P)\subset E$ by the $G$-action, defines a narrow element $P'\in E'(\ZZ)$ of order two
whose image in $E'(\FF_2)$ is non-zero. According to Theorem \ref{classification order two},
we may assume that $E'$ is given by a   Weierstra\ss{} equations of the form
$$
y^2+xy = x^3 + 4nx^2 +nx 
$$
for some sign  $n=\pm 1$, such that $P'\subset E'$ is given by the equations $x=y=0$. 
According to Theorem \ref{additional 2-division}, we have another element $Q'\in E'(\ZZ)$ of order two that
is disjoint from the zero-section, namely $Q'=(-4d, 2d)$ with the value $d=-n$. In turn, $P'+Q'$ is the third element of order two,
which must hit the zero-section, and   generates a subgroup scheme
$\mu_2\subset E'$. Since $R\subset E$ is disjoint from the zero-section, we infer
that $\mu_2=E[2]/G$.
In turn, we get an identification of elliptic curves
$E_\QQ=(E/E[2])_\QQ = (E'/\mu_2)_\QQ$. The latter quotient can be computed with Proposition \ref{quotient model},
and is    given by the Weierstra\ss{} equation
$y^2+xy=x^3+ a_2x^2 + x$
with coefficient $a_2=4n+6d=-2n$.  

Conversely, if $E_\QQ$ is given by one of the Weierstra\ss{} equation  $y^2+xy=x^3+ 2nx^2 + nx$ with $n=\pm 1$,
then conditions (i)--(ii) hold by Proposition \ref{equation sufficient}.
One easily sees that $R=(-n,n)$ has order four in in the groups $E(\QQ)$ and $E(\FF_2)$.
The statements on the singular fibers follow from Proposition \ref{fiber}.
Since all fibers are irreducible, every element in $E(\ZZ)$ is narrow.
The $j$-invariant is easily  computed with \cite{Magma}.
\qed

\section{Reinterpretation in terms of  stacks}
\mylabel{Reinterpretation}
 
In this final section we reformulate our main results in terms of stacks.
For details on the theory of stacks, we refer to the monographs of
Laumon and Moret-Bailly \cite{Laumon; Moret-Bailly 2000} and Olsson \cite{Olsson 2016}.

Let $g,r\geq 0$ be integers  with $g+r\geq 2$, and  $\shM_{g,r}$ be the 
Deligne--Mumford  stack over the base ring $\ZZ$ whose fiber categories
$\shM_{g,r}(R)$ are tuples $(C,\sigma_1,\ldots,\sigma_r)$ as follows:
$C$ is a finitely presented flat proper $R$-scheme
whose fibers are smooth curves with $h^0(\O_{C_a})=1$ and $h^1(\O_{C_a})=g$, and 
$\sigma_1,\ldots,\sigma_n$ are pairwise disjoint sections for the structure 
morphism $f:C\ra \Spec(R)$.
This has a natural compactification $\shMbar_{g,r}$, where the fibers $C_a$ have at most normal crossing singularities,
the sections pass through the smooth locus, and the automorphism group is finite.
The structure morphism $\shMbar_{g,r}\ra\Spec(\ZZ)$ is proper,
and the inclusion $\shM_{g,r}\subset\shMbar_{g,r}$ is open.
By abuse of notation, we write $\O_C(\sigma_i)$ for the invertible sheaf attached to the
image of the section $\sigma_i:\Spec(R)\ra C$, which is an effective Cartier divisor.
 
We are mainly interested in the case $g=1$ and $r\geq 1$.
Note that the objects $(C,\sigma_1)$ from $\shM_{1,1}$ can be regarded as   \emph{families of elliptic curves}
$E\ra\Spec(R)$, where the  zero-section $O\subset E$ is the image of $\sigma_1$.
Likewise, objects from $\shMbar_{1,1}$ can be regarded as families $Y\ra\Spec(R)$ of cubic
curves with a zero-section $O\subset Y$,   locally given by Weierstra\ss{} equations where $\Delta,c_4$
generate the unit-ideal. 
Given an object $(C,\sigma_1,\ldots,\sigma_r)\in\shMbar_{1,r}(R)$, the invertible sheaf $\O_C(\sigma_1)$
is semiample, and the homogeneous spectrum 
$$
P=P(C,\sigma_1)=\Proj H^0(C,\bigoplus\O_C(t\sigma_1))
$$
together with the induced section $\sigma_1:\Spec(R)\ra Y$ passing through $\Reg(Y/\ZZ)$
defines an object $(P,\sigma_1)$.
The construction $(C,\sigma_1,\ldots,\sigma_r)\mapsto (P,\sigma_1)$ actually yields a morphism of stacks
$\shMbar_{1,r}\ra\shMbar_{1,1}$, and composition with the $j$-invariant gives  
$\shMbar_{1,r} \stackrel{j}{\ra} \PP^1$.
By abuse of notation, we write $j(C)\in\PP^1(R)$ for the image of an object $(C,\sigma_1,\ldots,\sigma_r)$.
If $R=k$ is a field, this may be interpreted as a  number $j(C)\in k\cup\{\infty\}$.

Now let $m,n\in\ZZ$ be integers with $n$ odd and $\gcd(4m+1,n)=1$,
and consider our   Weierstra\ss{} equation $y^2+xy=x^3+mx^2+nx$. The Weierstra\ss{} model
$Y$ of the ensuing elliptic curve $E_\QQ$ together with its zero-section $O\subset Y$
can be regarded as an object in $\shMbar_{1,1}(\ZZ)$. This defines a map
$$
\Psi\lra \shMbar_{1,1}(\ZZ),\quad (m,n)\longmapsto (Y,O),
$$
defined on the set $\Psi=\{(m,n)\in\ZZ^2\mid \text{$n$ odd and $\gcd(4m+1,n)=1$}\}$.
The equations $x=y=0$ define  an element $P\in E(\ZZ)$ of order two. 
The strict transforms of $O\cup P\subset E$ on the minimal model $X$
yields a relatively semi-ample invertible sheaf $\shL$ on $X$.
The homogeneous spectrum   $Y'=P(X,\shL)=\Proj\bigoplus_{t\geq 0} H^0(X,\shL^{\otimes t})$
defines a partial resolution $X\ra Y'\ra Y$. Write $O',P'\subset Y'$ for the strict transforms of $O,P\subset Y$.
Each  geometric fiber  for the structure morphism $Y\ra\Spec(R)$ is either an   elliptic curve or a cycle of rational curves
with one or two irreducible components, and $P',Q'$ are contained in $\Reg(Y'/\ZZ)$.
This gives a map
$$
\Psi\lra \shMbar_{1,2}(\ZZ),\quad (m,n)\longmapsto (Y',O',P').
$$
If the element $P\in E(\ZZ)$ is narrow, which means $n=\pm 1$ according to Proposition \ref{narrow}, 
the partial desingularization coincides with the Weierstra\ss{} model, and the   map
becomes $(m,n)\mapsto (Y,O,P)$.

\begin{proposition}
The above map  induces a bijection between   $\{(m,n)\mid n=\pm 1\} $ 
and the set of isomorphism classes of objects $(C,\sigma_1,\sigma_2)\in\shMbar_{1,2}(\ZZ)$
with the properties
$\O_C(2\sigma_2)\simeq\O_C(2\sigma_1)$ and $j(C\otimes\FF_2)\neq \infty$.
\end{proposition}

\proof
Suppose first that $(C,\sigma_1,\sigma_2)$ arises from a pair $(m,n)$ with $n=\pm 1$.
Then $C=Y$, and the closed fiber at $p=2$ has invariant $j=1$.
Moreover, all geometric fibers of the structure morphism $C\ra\Spec(\ZZ)$ are integral.
Thus the Picard scheme $\Pic_{C/\ZZ}$ exists (\cite{FGA V}, Theorem 3.1), and we have an identification 
$E=\Pic^0_{C/\ZZ}$ of relative group schemes, given by $D\mapsto\O_C(D-O)$.
Since $P$ has order two in $E(\ZZ)$, the invertible sheaf $\O_C(2\sigma_1-2\sigma_2)$ has trivial
class in $\Pic_{C/\ZZ}(\ZZ)$. The Leray--Serre spectral sequence gives an exact sequence
$$
\Pic(\ZZ)\lra \Pic(C)\lra \Pic_{C/\ZZ}(\ZZ),
$$
and the factoriality of the ring $\ZZ$ ensures that $\O_C(2\sigma_1-2\sigma_2)$ is trivial.

Conversely, suppose that $(C,\sigma_1,\sigma_2)$ has the stated properties.
Then $j(C\otimes\FF_2)=1$, and this 
ensures that the generic fiber $E_\QQ=C_\QQ$ becomes an elliptic curve, with origin given by $\sigma_1$.
Let $E$ be its N\'eron model. Write $X$ and $Y$ for the minimal model and Weierstra\ss{} model, respectively.
From $\O_C(2\sigma_2)\simeq\O_C(2\sigma_1)$ we infer that $\sigma_1,\sigma_2:\Spec(R)\ra C$ pass through the same irreducible component,
in all fibers for $C\ra\Spec(\ZZ)$. It follows that 
the fibers are irreducible. Moreover, using that the sections pass through $\Reg(C/\ZZ)$
we infer   $C=Y$.  
 The image of $\sigma_2$ defines
a section $P\subset E$ that is disjoint from the zero-section $O\subset E$. Moreover, the element $P\in E(\ZZ)$ has order two
and is narrow. From Theorem \ref{equation necessary} and Proposition \ref{narrow} we see that $C=Y$ is given by a Weierstra\ss{} equation
$y^2+xy=x^3+mx^2\pm x$. The equation is unique  by Lemma \ref{unique equation}.
\qed
 
\medskip
Now consider the situation of Theorem \ref{classification order four}, such that we have a narrow element $R\in E(\ZZ)$
of order four whose image in $E(\FF_2) $ also has order four. The resulting narrow element $P=2R$
has order two in both $E(\QQ)$ and $E(\FF_2)$, and the three sections $O,R,P$ are pairwise disjoint.
This gives a map
$$
\{ (4,1),(-4,-1)\} \lra\shMbar_{1,3}(\ZZ),\quad (m,n)\longmapsto (Y,O, R, P).
$$

\begin{proposition}
The above map  induces a bijection between    $\{ (4,1),(-4,-1)\} $ 
and the set of isomorphism classes of objects $(C,\sigma_1,\sigma_2,\sigma_3)\in\shMbar_{1,3}(\ZZ)$
with the properties
$$
\O_C(4\sigma_2)\simeq\O_C(4\sigma_1),\quad \O_C(2\sigma_2)\simeq\O_C(\sigma_3)\quadand j(C\otimes\FF_2)\neq \infty.
$$
\end{proposition}

\proof
It follows from Theorem \ref{classification order four} that $m=4n$ and $n=\pm1$ gives a family 
of pointed stable curves $(C,\sigma_1,\sigma_2,\sigma_3)$ with the stated properties.
The two Weierstra\ss{} equations have different $j$-invariants, so the map in question is injective.

Conversely, suppose that  we have an object $(C,\sigma_1,\sigma_2,\sigma_3)$ with the properties at hand. 
Then $C$ and $\sigma_1$ define an elliptic curve $E_\QQ$.
Write $E$ for its N\'eron model, and let $R,P\in E(\ZZ)$ be the elements corresponding to $\sigma_2,\sigma_3$.
Set $S=\Spec(\ZZ)$. 
The conditions on the invertible sheaves $\O_C(\sigma_i)$ ensure $4R=O$ and $2R=P$. 
If a geometric fiber for $C\ra\Spec(\ZZ)$ contains a copy of $\PP^1$, it must intersect some image  $\sigma_i(S)$,
by stability. Again by the condition on the invertible sheaves $\O_C(\sigma_i)$ we infer that all geometric fibers
are irreducible. It follows that $C$ coincides with the Weierstra\ss{} model $Y$ of the elliptic curve $E_\QQ$.
This ensures that the structure morphism $E\ra S$ is semi-abelian.
Moreover, the subschemes $O,R,P\subset E$ map to the   $\sigma_i(S) \subset Y$ under the canonical morphism $E\ra Y$.
It follows that $O,R,P$ are pairwise disjoint. Hence $R$ has order four in both $E(\QQ)$ and $E(\FF_2)$.
Applying Theorem \ref{classification order four}, we see that the family of cubics $Y=C$ is define by the Weierstra\ss{} equation
$y^2+xy=x^3+4nx^2+nx$ with either $n=1$ or $n=-1$. 
\qed


\end{document}